\documentclass[11pt]{amsart}
\usepackage{amsfonts}
\usepackage{amsmath,amsthm}
\usepackage{latexsym}
\usepackage{amssymb}
\usepackage{enumerate}
\usepackage{comment}

\usepackage{hyperref}
\hypersetup{colorlinks,
	citecolor=blue,
	linkcolor=blue,
	urlcolor=blue}

\DeclareFontFamily{U}{wncy}{}
\DeclareFontShape{U}{wncy}{m}{n}{<->wncyr10}{}
\DeclareSymbolFont{mcy}{U}{wncy}{m}{n}
\DeclareMathSymbol{\Sh}{\mathord}{mcy}{"58}

\usepackage[notcite,notref]{showkeys}

\usepackage[english]{babel}
\usepackage{color}
\usepackage[latin1]{inputenc}

\numberwithin{equation}{section}
\usepackage{cite}

\theoremstyle{plain}
\newtheorem{theorem}{Theorem}[section]
\newtheorem*{theorem*}{Theorem}

\theoremstyle{remark}
\newtheorem{remark}[theorem]{Remark}

\newtheorem*{lem*}{Lemma}

\newtheorem*{remark*}{Remark}
\newtheorem*{NB*}{NB}

\newcommand{\R}{ \mathbb{R} }

\newcommand{\C}{ \mathbb{C} }
\newcommand{\Z}{ \mathbb{Z} }

\newcommand{\T}{ \mathbb{T} }

\newcommand{\vp}{\varphi}

\newcommand{\om}{ \omega }

\renewcommand{\phi}{ \varphi }

\newcommand{\eps}{\varepsilon}

\newcommand{\be}{\begin{equation}}
\newcommand{\ee}{\end{equation}}
\newcommand{\ben}{\begin{equation*}}
\newcommand{\een}{\end{equation*}}

\newcommand{\la}{ \langle }
\newcommand{\lan}{ \langle }
\newcommand{\lann}{ \langle\langle }
\newcommand{\llan}{ \langle\langle }
\newcommand{\ra}{ \rangle}
\newcommand{\ran}{ \rangle}
\newcommand{\rann}{ \rangle\rangle}
\newcommand{\rran}{ \rangle\rangle}
\newcommand{\p}{ \partial}

\title[Stochastic perturbations of  integrable systems]
{ Two approaches to average  stochastic perturbations of  integrable systems. }

\author{Sergei Kuksin}
\address{ Universit\'e Paris Cit\'e and Sorbonne Universit\'e, CNRS, IMJ-PRG, F-75013 Paris, France 
 \&  Steklov Mathematical Institute of RAS, 
119991 Moscow, Russia \&
 Peoples' Friendship University of Russia (RUDN University),
117198, Russia  \& National Research University Higher School of
	Economics, 119048 Moscow, Russia}
\email{ sergei.kuksin@imj-prg.fr}

\begin{document}
	
\begin{abstract}
We discuss two approaches to study the long-time behaviour and infinite-time behaviour 
of solutions for integrable hamiltonian systems under small stochastic perturbations. Then
we compare these results with those for deterministic perturbations of integrable systems. 
\end{abstract}

\maketitle

  \section{The setting}

 In symplectic space $(\R^{2n}_{x,y}, dx\wedge dy)$ we consider a Birkhoff--integrable Hamiltonian system --  its 
 Hamiltonian $H$ depends only on the actions   $I_j = (x_j^2 +y_j^2)/2$. Introducing complex coordinates 
 $v_j = x_j +iy_j$  we write this  system in the convenient complex form:
 \be\label{Birk}
 \tfrac{\p}{\p t} v_k(t) = i \nabla_k H(I) v_k, \;\;\; k=1,\dots, n, \quad I=I(v), \quad v(t) \in \C^n.
    \ee
     Note that $I_j = \frac12 |v_j|^2$, and that 
      in $\C^n$ we use the real scalar 
   product 
   $
   (v, w) \mapsto
   \Re \la v, \bar w\ra = \Re \sum v_j \bar w_j. 
   $
   Study of small perturbations of this system for large values of time, either globally in $\R^n$, or locally near an 
   equilibrium, is a classical problem of dynamical systems starting  the end on XVIII century.\footnote{Another 
   classical question is to study solutions of perturbed systems near the tori $\{v: I_k(v)=I_k^0\ne0, k=1, \dots, n\}$,
   which are invariant for system \eqref{Birk}. We do not touch it here.   }
   There are two main settings:
   
   A) the perturbation is a small smooth vector field,
   
   B) the perturbation is a small Hamiltonian field.
   
   We will consider a third case and then compare the obtained results with what is known for the cases  A) and  B).

Consider  an $\eps$-small stochastic perturbation of the integrable  system \eqref{Birk}. Firstly, only 
 for $0\le t\le T\eps^{-1}$:
 \be\label{1} 
  \tfrac{\p}{\p t}  v_k( t) = i  {\nabla_k H} (I) v_k+  \eps^{}P_k(v) +\sqrt{\eps} \sum_{j=1}^{n_1} B_{k j}(v) 
    \dot \beta_j^c( t), 
 \;\;\;  \qquad v(0)= v^0,
   \ee
   where $k=1,\dots, n$, 
    and  $\{ \beta_k^c( t), 1\le k\le n\}$ are standard independent complex Wiener processes. That is, 
   $\beta_k^c( t)=\beta_k^+( t)+i \beta_k^-( t)$, and  $\beta_k^\pm(t)$ are standard  independent  real 
    Wiener processes.   It is convenient to
pass to the slow time $\tau =\eps t$ and re-write the  $v$-equations above  as
 \be\label{2} 
 \dot  v_k(\tau) = i \eps^{-1} \nabla_k H (I) v_k+ P_k(v) + \sum_{j=1}^{n_1} B_{k j}(v)   \dot \beta_j^c( \tau), 
  %b_k \dot\beta_k^c(\tau), 
 \;\;\; k=1,\dots, n,    
 \ee
 $v(0)= v^0$. Here and below the upper dot stands for $d/d\tau$, and $\{\beta_j^c(\tau)\}$ is another 
 collection of  standard independent complex Wiener processes. 
  Or, in vector notation,  as 
  \be\label{v_eq} 
 \dot  v(\tau) = i \eps^{-1} \text{diag}\{ {\nabla_k H (I) }\} v+ P(v) + B(v) \dot\beta^c(\tau),
%  \{b_k \dot\beta_k^c(\tau)\},
\quad  v(0)= v^0,
  \ee
  where  $ v(\tau) \in \C^n$ and %$0\le \tau\le  T$. 
  $0\le \tau\le T$. A solution of this equation will be denoted $v^\eps(\tau)$. 
  
   We impose the following assumptions on equation \eqref{v_eq}:
   
\noindent  {\bf (A)} 1) (Anosov's non-degeneracy):
For a.a.  $I\in\R^n$  components of  the frequency  vector $\nabla H(I)$ are
independent over integers. That is, 
 \be\label{Anosov} 
\text{if\  $s\in \Z^n$ \ and\  $\nabla H(I) \cdot s=0$, \; then\  $s=0$.}
\ee
2) Vector-functions  $\nabla H(I)$, $P(v)$ and operator-function 
$B(v)$  are locally Lipschitz, with at most polynomial growth at infinity.
\\
 3) The $v$-equation \eqref{v_eq}    is well posed for $0\le\tau\le T$ and there the first few moments of 
its solutions are bounded uniformly in $\eps$.\\
 4) The noice in the equation 
  is non-degenerate: the rank of the complex matrix  $B(v)$ equals $n$ for all $v$ (so $n_1 \ge n$).

\begin{remark}
  Assumption 1) (Anosov's non-degeneracy) holds  e.g. if $\nabla H$ is a constant vector-function with some 
   components $\lambda_j$, independent over integers.  Then the  unperturbed system \eqref{Birk} 
   is  linear,
  \be\label{linear} 
 \tfrac{\p}{\p t}  v_k = i \lambda_k v_k, \quad k=1, \dots, n. 
 \ee
 The theorems below apply to this case, 
  but their stronger versions may be found in our previous work  \cite{HK}.  The results of   \cite{HK}
  apply when the matrix $B(v)$ in \eqref{v_eq} is degenerate, and their suitable versions hold if the 
  frequencies $\{\lambda_j\}$ are dependent over integers.  
It is easy to see that  1) also holds  if $H$ is $C^2$-smooth and 
  $$
  \text{meas} \{I: \det d^2H (I) =0\} =0
  $$
  (this is the Kolmogorov nondegeneracy condition). 
  
  Assumption  3) (the apriori estimates)  is valid, for example, 
  if the matrix $B(v)$ is uniformly bounded,  and   
 $$
 \Re \lan  P(v), \bar v\ran  \le  C |v|^2 + C_1 \qquad \forall v\in \C^n.
 $$  
 \end{remark}

      \section{ First approach:
      the fast-slow system.       }

 Let us introduce in $\C^n=\{v=(v_1, \dots, v_n)\}$  the usual   action-angle coordinates $(I,\vp)$, where 
  $I=(I_1, \dots, I_n) \in \R_+^n$ and  $\vp= ( \vp_1, \dots, \vp_n)  \in \T^n = \R^n/(2\pi\Z)^n$: 
 $$
 I_j = \tfrac12 |v_j|^2,  \quad \vp_j = \arg v_j, \;\;\; j=1,\dots, n
    $$ 
    (so $v_j = \sqrt{2I_j} \exp{i\vp_j}$). 
    Then   $v$-equation \eqref{v_eq} 
     may be re-written  in terms of these coordinates. The equations for the actions are:
      \be\label{3} 
 \dot  I^{\eps}_k(\tau) = \Re ( \bar v^{\eps}_k\, P_k(v^{\eps}))  +  \sum_l |B_{kl}(v^\eps)|^2
  + \sum_{j=1}^{n_1}  \Re\big( \bar v^{\eps}_k B_{kj}(v^{\eps})\big)
 \dot\beta_j(\tau),
   \ee
where $k=1,\dots, n$ and 
$\beta_1(\tau), \dots, \beta_{n_1}(\tau)$  are independent standard real  Wiener processes. (Recall that 
 in $\C$ we use the real scalar product $(z_1,z_2) \mapsto \Re( z_1, \bar z_2)$). 
Re-denoting  the terms of these equations  and using vector notation, we  write   system \eqref{3}
 in the vector form:
 \be\label{I_eq} 
 \dot  I^{\eps}(\tau) = F( v^{\eps}) + G(v^\eps) \dot\beta(\tau), \qquad I^\eps(0) = I^0 := I(v^0),
 %  + \sum_{j=1}^{n_1}  \Re( \bar v^{\eps}_k B_{kl}(v^{\eps}))
   \ee
   where $0\le \tau \le T$. 
The $\vp$-equations are:
 $$ %\be\tag{$\vp$-equations}
\dot \vp^{\eps}_k (\tau) = \eps^{-1} {\nabla_k H}(I^{\eps}) + \text{  term of order one with   strong singularity at} \;
 \{I^{\eps}_k=0\},
$$
$1\le k\le n$. 
We have got a  {\it slow-fast system}, where the actions 
$I^{\eps}_k$'s are {\it slow variables} and the angles  $\vp^{\eps}_k$'s are {\it fast variables. } As usual, we care the most
about the actions. 
The logic of averaging suggests that their limiting  behaviour  as $\eps\to0$ 
can be described by  averaging  in  angles $\vp$ of the  slow system \eqref{I_eq}.

Accordingly, let  us  average \eqref{I_eq} in the fast variables $\vp$, 
 using the rules of  stochastic  averaging:  
   \be\label{I_aver} 
 \dot  I(\tau) =\lann   F \rann(I) + \lann G\rann (I)  \dot\beta(\tau),
 \;\;\;  I(0) = I^0, \quad I(\tau) \in \R_+^n.
   \ee
   Here  $\beta(\tau)$ is  the standard   Wiener process  in $\R^n$.  The vector field
      $\lann F \rann(I)$  is the averaging in    $\vp\in \T^n$ of the  field  $F$, written as       $F(I,\vp)$. I.e.,
       \be\label{F}
   \lann F \rann(I) = \tfrac1{(2\pi)^n} \int_{\T^n}   F(I,\vp)\,d\vp.
  \ee
       Since 
  $\ 
   v_j =r_j e^{i\vp_j}$, 
   where 
   $r_j =\sqrt{2I_j}$,
  and $F$ is locally Lipschitz in $v$, then $ \lann F \rann$ is a locally Lipschitz  vector-function of $r=(r_1, \dots, r_n)$. So this is  a
  H\"older-$\frac12$ function of $I$. 
  The  real $n\times n$-matrix  $ \lann G\rann (I)$ is an averaging of $G(I,\vp)$,  defined in a related, but 
   more complicated way, which  still is explicit;   see \cite{HKP}  and references there. A-priori it also is only  H\"older-$\frac12$ smooth. But 
  \be\label{const_diff}
  \begin{split}
   \text{if $B$ is a constant matrix, then $ \lann G \rann =\,$diag$\, \big(b_k\sqrt{2I_k}\big)$, }\\
   \text{
     where all $b_k$ are positive numbers}.
     \end{split}\ee
     
       A partial description of the limiting as $\eps\to0$
        dynamics of actions in eq.~\eqref{v_eq}  is
     given by the following result (for a proof see \cite[Section 6]{HKP}, and see \cite{FW} for a similar assertion). 
      
     \begin{theorem} \label{t_1}
      The collection of  curves of  actions $I^\eps(\tau)\  (0\le\tau \le T)$
      of  solutions  $v^\eps(\tau)$ for \eqref{v_eq} 
       with $0<\eps\le1$ is    precompact in the space of curves, in the 
     sense of distribution. 
     
     1) Every limiting (in the sense of distribution) random process 
     $I^0(\tau)$ is a weak solution of the averaged equation \eqref{I_aver}.
      
     2) Eq. \eqref{I_aver}  always has at least one 
      solution $I^0(\tau)$. If such a solution is unique,   then  $I^\eps(\tau)$ converges 
     to $I^0(\tau)$ in the sense of distribution as $\eps\to0$. 
      \end{theorem}

     We recall that in \eqref{I_aver} $\lann   F \rann$ and $\lann   G \rann$ are H\"older-$\frac12$ function of $I$, 
     so there is no reason for this equation 
       to have only one solution. Related to that we have the following facts  (see \cite{HKP}):

      i) A result of H.Whitney \cite{HW} implies that 
      if the coefficients of   the $v$-equation  are $C^2$-smooth, then $\lann   F \rann(I)$ is $C^1$-smooth, 
      and 
      $ \lann   G \rann(I) = \sqrt{K(I)}$, where matrix  $K(I)$  is $C^1$-smooth, selfadjoint non-negative, 
      and \ det$\,K(I)=0$ if some $I_j$ vanishes.

      \noindent
     This is better than before, but still $ \lann   G \rann(I)$ has a $\sqrt{I_j}$-singularity when $I_j=0$, and it is not 
     clear if  solutions for equations like that are unique.

     ii) Let in  \eqref{v_eq}  the dispersion matrix  $B$ is constant. Then by \eqref{const_diff}
     the noise in  \eqref{I_aver} 
     is given by a (non-constant) 
     diagonal matrix which has       square-root singularities when some $I_j$ vanishes. 
      For this case the result of the paper   \cite{YW} 
     implies that eq.  \eqref{I_aver}  has a unique solution. So item 2) of Theorem~\ref{t_1}   applies, and  \eqref{I_aver}  describes 
     the limiting behaviour of actions of solutions $v^\eps(\tau)$ 
      for $0\le \tau \le T$.
     
     But what happens  if the  matrix  $B(v)$ is not constant? Then  $ \lann   G \rann(I) = \sqrt{K(I)}$ 
      has  complicated H\"older-$\frac12$ singularities  when $I_j=0$ for some $j$.
      No uniqueness result is known for this case.  On the contrary, 
     
     iii) In paper  \cite{WY}     the authors give an example of a system with H\"older-$\frac12$ singularity at $I_j=0$ where a solution 
    is not unique.

\begin{remark}
In paper \cite{Cher} is considered a class of stochastic equations on a unit cube $K^n \subset \R^n$ with suitable 
boundary conditions at $\partial K^n$ and with the stochastic term $\sqrt {K(v)}\,  d\beta(t)$, where $K(v)$ is a smooth 
non-negative symmetric matrix which degenerates for some $v$. Using the techniques of martingale solutions it is established there 
that the equation has a unique weak solution. We do not know if the results of \cite{Cher} may be extended to stochastic 
equations in $\R^n$ and to equations in $\R^n_+$ of the form \eqref{I_aver}, and are not aware of publications, related to this 
question. 
\end{remark}

    {\bf   Conclusion.}  Equation  \eqref{I_aver}  describes the limiting dynamics of 
      actions of solutions $v^\eps(\tau)$  as $\eps\to0$, only if in \eqref{v_eq} 
      $P(v)$ is $C^2$-smooth and $B$   is a constant matrix. 
      \smallskip
      
      What can we do in the general case?

      \section{  Second approach: the effective equation.}\label{s_3}

 Let us come back to the original  $v$-equation  \eqref{v_eq}.
   I recall that the celebrated Krylov--Bogolyubov averaging for  perturbations of linear conservative systems in 
  $\R^{2n}$ (see in \cite{AKN}) 
   is performed in $\R^{2n}$ WITHOUT passing to the action-angles. It describes the limiting dynamics 
  of actions via
  certain  {\bf effective equation} in $\R^{2n}$, derived using the {\bf interaction representation}.\footnote{Which is 
   some non-autonomous linear change of variables in $\R^{2n}$; e.g. 
  see in \cite{HK, GK}.}
  
   The latter does not   exist for our $v$-equation \eqref{v_eq}, where the unperturbed equation \eqref{Birk} is
   non linear. 
   Still, motivated by  Krylov--Bogolyubov we GUESS a suitable effective equation. 
    It is  obtained from the $v$-equation by:
    
-- dropping its first term, 

-- suitably averaging the other two terms. \\
So the equation is 
 \be \label{eff_eq}
 \dot a(\tau) = \llan P\rran (a)  + \llan B\rran (a) \dot\beta^c(\tau),
\qquad a(\tau) \in \C^n, \quad
 a(0)\!=\! v^0. 
  \ee
Here the averaging of the vector field $P$  is 
  $$
   \lann P \rann(a) = \tfrac1{(2\pi)^n} \int_{\T^n} \Phi_\om\circ  P(\Phi_{-\om} a)\,d\om,
  $$
where $\Phi_\om$ is the rotation operator in $\C^n$, \ 
  $
   \Phi_\om = \text{diag} \big(e^{i \om_1}, \dots, e^{i\om_n} \big), \; \om\in \T^n.
  $
  
   The averaging $\llan B\rran (a)$ of the matrix-field $B(a)$ is more complicated, but still is explicit. 
   This $\llan B\rran$   is a non-degenerate   complex $n\times n$-matrix.  
  Vector field $ \llan P\rran(a)$ and matrix field $\llan B\rran(a) $ are as smooth as the components 
   $P(v)$ and $B(v)$ of the $v$-equation.  See in \cite{HKP}. 
 
 \begin{theorem} \label{t_2}
    Under the imposed assumptions A1)--A4)    effective equation \eqref{eff_eq} 
      has a unique solution $a(\tau)$, which inherits the estimates in A3). 
      When $\eps\to0$,  the actions $I^\eps(\tau)= 
   I(v^\eps(\tau))$ of solutions for  $v$-equation \eqref{v_eq} 
    with  $0\le\tau\le T$ converge in distribution to the action-vector   $I(a(\tau))$. The curve $I(a(\tau))$ satisfies equation 
    \eqref{I_aver}. 
   \end{theorem}

Now in addition to  assumptions  A1)--A4) let us assume that:
\medskip

\noindent  A5)  $v$-equation \eqref{v_eq} is well posed for $\tau\ge0$, and  the first few moments 
of its solutions are bounded uniformly in $\eps$ and in $\tau$ 
(this is an amplification of assumption  A3)\,).

  \begin{remark}\label{r_2}
   Assumption A5) holds if $B(v)$ is uniformly bounded, and $P(v)$ is coercive:  
 \be\label{coerc}
 \Re \lan P(v), \bar v\ran \le -\alpha_1 |v| + \alpha_2, \quad \alpha_1>0. 
 \ee
(This means that the vector field $P$ ``contains friction".)
\end{remark}

 \begin{theorem} \label{t_3}
   If assumptions  A1)-A5) are satisfied,     then the solution $a(\tau)$ of effective equation \eqref{eff_eq}
    meets  the estimates in A3) uniformly in $\tau\ge0$, 
   and the   convergence in Theorem \ref{t_2} also holds uniformly in time $\tau$. 
   \end{theorem}
   
   The proof of this result crucially uses the fact that since the first moments of solutions $a(\tau)$ are 
   bounded uniformly in time and matrix $\llan B\rran (a)$ is non-degenerate, then equation \eqref{eff_eq} is 
   mixing: in $\C^n$ exists a measure such that for any $v^0$ the distribution of solution $a(\tau)$ 
   converges to it, as $\tau\to\infty$.   See \cite[Section~10]{HKP} and \cite[Section 7]{HK}.

   I recall that a vector field $\tilde P$ in $\C^n$ is hamiltonian with a $C^1$-smooth   {\bf real} Hamiltonian $h$, 
   if 
 $$
 \tilde P_k(v) = 2i \tfrac{\p}{ \p \bar v_k}  h(v), \qquad k=1, \dots, n. 
 $$
 Assume that in   \eqref{v_eq} 
 $\;
P(v) = P_1(v) + P_2(v),
 $ 
 where the field $P_2$ is hamiltonian. 
 
  \begin{theorem} \label{t_4}
  If $P(v)$ is as above, then the effective equation may be modified by removing 
   the averaging of the    hamiltonian part $P_2$: 
 \be \label{eff_eq_modif}
 \dot a(\tau) = \llan P_1\rran (a)  + \llan B\rran (a) \dot\beta^c(\tau),
\qquad  v^{\eps}(0)\!=\! v^0. 
  \ee
  That is, the  assertions of Theorems \ref{t_2} and \ref{t_3}  remain true if for $a(\tau)$ one takes  a solution of equation
 \eqref{eff_eq_modif}. 
\end{theorem}

In \cite{GK} this result is established for perturbations of linear systems \eqref{linear}. A proof for solutions of
\eqref{v_eq} is similar and uses in addition the construction from \cite[Section~8]{HKP}.

   \section{Local versions of the results. }\label{s_5}
  Our starting assumption was that the unperturbed integrable system can be  put to the global Birkhoff normal form
 \eqref{Birk}. This  is restrictive. Now let us assume that the Birkhoff normal 
 form  holds only locally, for $v\in Q$, where $Q$ is a neighbourhood of the origin in $ \C^n$ of the form 
 $\{v: |v_j|\ \le C_j>0, j=1, \dots, n\}$.\footnote{This may be achieved by using the Vey theorem, see \cite{Ell}.}
 Consider the same perturbed $v$-equation as before:
 $$
 \dot  v(\tau) = i \eps^{-1} \text{diag}\{ {\nabla_k H (I) }\} v+ P(v) + B(v) \dot\beta^c(\tau),
%  \{b_k \dot\beta_k^c(\tau)\},
\quad  v(0)= v^0 \in Q.
  $$
  Then the assertions of  Theorems \ref{t_2} and \ref{t_4}  hold for $0\le \tau \le \theta$, where $\theta$ is the  random   moment of time 
  when a solution of 
    \eqref{v_eq}  exits $Q$. This is a random variable of order one:
   $$
   \text{probability that } \ \theta\le\lambda  \quad  \text{is bounded by }\;
   C\sqrt\lambda, \quad \forall\,\lambda>0. 
      $$
See \cite[Section 9]{HKP}.

   \section{ Applications.     }
   
   A more detailed discussion of applications of out results to the equations below see in \cite{HKP}. 
 
\subsection{Chains  of oscillators}  For the unperturbed  integrable  system \eqref{Birk} let us  
take a   chain of  nonlinear oscillators:
 \be\label{chain}
 \ddot q_k = -Q(q_k), \quad k=1, \dots, n, 
   \ee 
   where the potential  $Q$ is a smooth  function, which is convex and  odd.  Assume that it is  ``nice" at infinity, while  near zero
   $$
   Q(q) = \alpha q + \beta q^3 +o(|q|^4),   \quad \alpha, \beta>0.
   $$
   Then this  system is globally integrable and may be written in the form \eqref{Birk}. 
   Theorems \ref{t_2}-\ref{t_4}  apply to study its stochastic perturbations. The results obtained in this way are 
   related to the non-equilibrium statistical physics.

\subsection{Damped/driven hamiltonian systems}  Consider the following $v$-equation in $\C^n$:
 \be\label{11}
 \dot v_k(\tau)  = i\eps^{-1} {\nabla_k H}\big(I(v)\big) v_k +2i \tfrac{\p}{ \p \bar v_k}  h(v) - \nu_k v_k + b_k \dot\beta_k^c(\tau), \quad
 1\le k\le n,
   \ee 
   where $\nu_k>0$ and $b_k\ne0$ for all $k$. Let $H, h \in C^2$, the frequency map $I\mapsto \nabla H(I)$ is
   Anosov-nondegenerate and A5) is fulfilled.\footnote{In particular, A5) holds if the hamiltonian field
   $i \nabla_{\bar v} h(v)$ meets \eqref{coerc} with $\alpha_1=0$.}
   Then the assumptions  A1)-A5) hold. The averaging does not change the vector field
   $
   \dot v = -$diag$\{\nu_k\} v,
   $
   so in  view of
   Theorem~\ref{t_3},   for a modified  effective equation one may take the system 
 \be\label{12}
 \dot a_k(\tau)  =- \nu_k a_k + b_k \dot\beta_k^c(\tau), \quad
 1\le k\le n.
   \ee 
   For small $\eps$ the actions of solutions for eq. \eqref{12} approximate (in distribution) 
    the actions of solutions for eq. \eqref{11}, uniformly in  time. Note that  eq. \eqref{12} may be integrated  explicitly. (Its solution is a collection of $n$
    independent complex Ornstein--Uhlenbeck processes.)

   \section{ Discussion}\label{s_6} First we  compare the   results above  with those, available for  non-random perturbations of integrable 
   Hamiltonian systems.

   \subsection{Hamiltonian perturbations of  integrable systems}  Consider the following deterministically
   perturbed  integrable   system, written in the slow time $\tau$:
   \be\label{15}
 \dot v_k(\tau)  = i\eps^{-1} {\nabla_k H}\big(I(v)\big) v_k +2i \tfrac{\p}{ \p \bar v_k}  h(v), \quad
 1\le k\le n, \qquad v(0) = v^0.
   \ee 
   If $H$ and $h$ are analytic, then    under the
     steepness assumption on Hamiltonian $H$, Nekhoroshev's  theorem (see in \cite{AKN}) 
     applies, provided that $v^0$ lies
      not too  close to the locus 
 $$
 X= \{ v\in \C^n: I_k(v) =0\;\; \text{for some}\;\; k\}.
   $$ 
   The theorem claims that  for times \ $\tau \lesssim \!\exp{\eps^{-a}}$, $a>0$,  actions of  solutions 
    for eq.~\eqref{15}  stay  $\eps^b$-close $(b>0)$  to    the original action-vector $I(v^0)$.    It means that 
    during exponentially long time the hamiltonian term of eq.~\eqref{15} does not affect the dynamics of 
    actions of its solutions, within  precision $O(\eps^b)$. 
   
    Theorem~\ref{t_3} above deals with perturbations \eqref{v_eq}  of  the integrable part of    eq.~\eqref{15},
   and for any initial data $v^0$ it
   allows to control actions of  solutions for  \eqref{v_eq}  uniformly in time. This happens since the perturbation in
    \eqref{v_eq} contains noice and friction (see Remark~\ref{r_2}) which stabilise the dynamics. In particular, if we add to each
    equation in ~\eqref{15} dissipation $-\nu_k a_k$ and randomness $b_k \dot\beta_k(\tau)$, then it 
    turns out that the hamiltonian term does not at all affect the dynamics of actions, within precision $o(1)$.

   We   note that in  Theorem~\ref{t_3} we had to  assume that the unperturbed system is
   integrable in the whole space, while for the validity of Nekhoroshev's  result  it suffices to assume its 
   integrability only in  some open domain which contains $v^0$ (and is diffeomorphic to $\T^n\times O$, where
   $O$ is a domain in $\R^n$). 
   
   Nekhoroshev's theorem does not allow to take initial data $v^0$ close to $X$ since its proof uses the 
   action-angle coordinates, which are singular at $X$, while Theorem~\ref{t_3} does not need them.  For some class of perturbed systems \eqref{15} Niederman in  \cite{LN} established a version of 
   Nekhoroshev's theorem in the vicinity of the origin, avoiding the action-angles.  So his result applies to
   solutions of \eqref{15} with small $v^0$. It seems that for the moment of writing there is no version of 
   Nekhoroshev's theorem, valid for  solutions with initial data  in $X$, or  close to it.

     \subsection{ Arbitrary perturbations of integrable systems
     } Now consider perturbed system
   \be\label{16}
 \dot v_k(\tau)  = i\eps^{-1} {\nabla_k H}\big(I(v)\big) v_k + P_k(v), \quad
 1\le k\le n, \qquad v(0) = v^0,
   \ee 
where $H$ and $P$ are $C^\infty$-smooth. Let $H$ satisfies 
   Bakhtin's nondegeneracy condition (see in \cite{AKN})  and  $v^0$ be a random vector, $v^0 = v^0(\om)$. 
    Then by the 
   Neishtadt--Bakhtin theorem for  random    parameters $\om$ 
     outside some set of small measure $\ \lesssim \eps^b$, $b>0$,\footnote{In particular,  the initial data 
     $v^0(\omega)$      should avoid the vicinity of $X$.}
      for $0\le \tau\le T$ actions of solutions for 
   eq.~\eqref{16} stay close to solutions for the averaged equation for actions:
      $$
 \dot  I(\tau) =  \lann   P \rann (I), 
 \qquad  I(0) = I^0,
   $$
   where $\lann   P \rann$ is defined as in 
    \eqref{F}.  This result is obviously related to our Theorem~\ref{t_2},  to its
   local version, discussed in Section~\ref{s_5}, and to Theorem~\ref{t_3}. 
    The assertion of  Theorem~\ref{t_3} is significantly stronger 
   than the result of Neishtadt--Bakhtin. 
    Again, it happens since the perturbed equation \eqref{v_eq}, examined in Theorem~\ref{t_3}, is
   stabilised by noise and friction (and since the theorem assumes the global Birkhoff coordinates). 
   
   \subsection{Related results} Versions of the results in Sections 3,\,4 are available for some classes 
   of stochastic perturbations of linear and  integrable  PDEs. See references in \cite{HK, HKP}. The  action-angle approach as in 
   Section~2, being applied for SPDEs, may allow to get for them the first assertion of Theorem~\ref{t_1}. But to establish
    that the corresponding  equation \eqref{I_aver} has a unique solution 
    (as is required for the validity of the second assertion of the theorem)    seems to be impossible. 
    
    If the unperturbed integrable system is the linear system  \eqref{linear}, then the results of Sections~3,\,4 are applicable to its
     stochastic perturbations if the frequencies $\{\lambda_k\}$ are independent over integers. But  suitable versions of those 
     results apply if the frequencies are any $n$ non-zero real numbers and the matrix $B(v)$ in \eqref{v_eq} is a locally Lipschitz 
     complex $n\times n_1$-matrix, see in \cite{GK}. 
     The action-angle approach seems to be non applicable in this general setting.

 \end{document}